\documentclass{amsart}

 \usepackage{amssymb,latexsym,amsmath,amsthm}

\usepackage{url}

\usepackage[numbers,sort&compress]{natbib}

\usepackage[numbers,sort&compress]{natbib}

\numberwithin{equation}{section}

\newtheorem{thm}{Theorem}[section]
\newtheorem{prop}{Proposition}[section]

\newtheorem{lemma}{Lemma}[section]

\newtheorem{remark}{Remark}[section]

\newtheorem{cor}{Corollary}[section]

\begin{document}

\author{Mostafa Fazly}
\author{Juncheng Wei}
\author{Xingwang  Xu}

\address{Department of Mathematical and Statistical Sciences, CAB 632, University of Alberta, Edmonton, Alberta, Canada T6G 2G1}

\email{fazly@ualberta.ca}

\address{Department of Mathematics, University of British Columbia, Vancouver, B.C. Canada V6T 1Z2.}
\email{jcwei@math.ubc.ca}

\address{Department of Mathematics, National University of Singapore, Singapore 119076, Singapore}
\email{matxuxw@nus.edu.sg}

\thanks{The first two authors are supported by NSERC grants.}

\title{A pointwise inequality for the fourth order Lane-Emden equation}

\keywords{Semilinear elliptic equations, a priori pointwise estimate, Moser iteration type arguments, elliptic regularity}

\maketitle

\begin{abstract}  We prove that the following pointwise inequality holds
 \begin{equation*}
-\Delta u \ge \sqrt\frac{2}{(p+1)-c_n} |x|^{\frac{a}{2}} u^{\frac{p+1}{2}} + \frac{2}{n-4} \frac{|\nabla u|^2}{u} \ \ \text{in}\ \ \mathbb{R}^n
 \end{equation*}
 where   $c_n:=\frac{8}{n(n-4)}$,  for positive bounded solutions of the fourth order H\'{e}non equation that is
\begin{equation*}
\Delta^2 u = |x|^a u^p  \ \ \ \  \text {in }\ \ \mathbb{R}^n
\end{equation*}
for some $a\ge0$ and $p>1$. Motivated by the Moser's proof of the Harnack's inequality as well as Moser iteration type arguments
in the regularity theory, we develop an iteration argument to prove
the above pointwise inequality. As far as we know this is the first
time that such an argument is applied towards constructing pointwise
inequalities for partial differential equations. An interesting
point is that the coefficient $\frac{2}{n-4}$ also appears in
the fourth order $Q$-curvature and the Paneitz operator. This in
particular implies that the scalar curvature of the conformal metric
with conformal factor $u^\frac{4}{n-4}$ is positive.
 \end{abstract}

\section{Introduction}

We are interested in proving a priori pointwise estimate for
positive solutions of the following fourth order H\'{e}non equation
\begin{equation}\label{4henon}
\Delta^2 u = |x|^a u^p  \ \ \ \  \text {in }\ \ \mathbb{R}^n
\end{equation}
where $p>1$ and $a\ge 0$.   Let us first mention that for the case
$a=0$, it is known that (\ref{4henon}) only admits $u=0$ as a
nonnegative solution when $p$ is a subcritical exponent that is
$1<p<\frac{n+4}{n-4}$ when $n\ge 5$ and $1<p$ when $n\le 4$.
Moreover,  for the critical case $p=\frac{n+4}{n-4}$ all entire
positive solutions are classified. See \cite{l,wx}.  This is a
counterpart of the standard Liouville theorem of Gidas-Spruck in
\cite{gs,gs2} for the second order Lane-Emden equation
\begin{equation}\label{2lane}
-\Delta u =  u^p  \ \ \ \  \text {in }\ \ \mathbb{R}^n
\end{equation}
stating that  $u=0$ is the only nonnegative solution for
(\ref{2lane}) when $p$ is a subcritical exponent that is
$1<p<\frac{n+2}{n-2}$ when $n\ge 3$.      Note also that for the
fourth order  H\'{e}non equation, it is conjectured that $u=0$ is
the only nonnegative solution of (\ref{4henon}) when $p$ is a
subcritical exponent that is when $1<p<\frac{n+4+2a}{n-4}$ and $n\ge
5$, see \cite{fg}.  Therefore, throughout this note, when we are
dealing with (\ref{4henon}), we assume that $p > \frac{n+4+2a}{n-4}$
and $n\ge 5$.  For more information, see \cite{fg,so} and references
therein.

Pointwise estimates have had tremendous impact on the theory of
elliptic partial differential equations.  In what follows we list
some of the celebrated pointwise inequalities for certain semilinear
elliptic equations and systems. These inequalities have been used to
tackle well-known conjectures and open problems.  The following
inequality by Modica \cite{m}  has been one of the main techniques
to solve the De Giorgi's conjecture (1978) for the Allen-Cahn
equation and to analyze various semilinear equations and problems.

\begin{thm} (Modica \cite{m}, 1985) Let $F\in C^2(\mathbb {R})$ be a nonnegative function and
$u$ be a bounded entire solution of
\begin{equation}\label{mod}
 \Delta u=F'(u) \ \ \text{in} \ \  \mathbb{R}^n. 
 \end{equation}
  Then
  \begin{equation}\label{pointmod}
  |\nabla u|^2\le 2 F(u) \ \ \text{in}\ \  \mathbb{R}^n.
  \end{equation}
\end{thm}
For the specific case $F(u)=\frac{1}{4}(1-u^2)^2$, equation (\ref{mod}) is known as the Allen-Cahn equation.
 Note also that Caffarelli et al.  in \cite{cgs} extended this
inequality to quasilinear equations.   We refer interested readers to \cite{fv1,fv2,fv3,fv4,cfv,fsv} regarding pointwise gradient estimates and certain improvements of (\ref{pointmod}).    For the fourth order counterpart of (\ref{mod}) with an arbitrary
nonlinearity, a general inequality of the form (\ref{pointmod}) is
not known.   However, for a particular nonlinearity known as the
fourth order Lane-Emden equation, i.e.
\begin{equation}\label{4lane}
\Delta^2 u =  u^p  \ \ \ \  \text {in }\ \ \mathbb{R}^n
\end{equation}
 it is shown by Wei and Xu, as Theorem 3.1 in  \cite{wx}, that the negative Laplacian of the positive solutions is non-negative
 that is $-\Delta u\ge 0$ in $\mathbb{R}^n$.   Set $v=-\Delta u$ and from the fact that  $-\Delta u\ge 0$ we can consider \eqref{4lane} as a special case (when $q=1$) of the Lane-Emden system that is
\begin{eqnarray}\label{lane-emden}
 \left\{ \begin{array}{lcl}
\hfill -\Delta u&=& v^q   \ \ \text{in}\ \ \mathbb{R}^n,\\
\hfill -\Delta v&=& u^p   \ \ \text{in}\ \ \mathbb{R}^n,
\end{array}\right.
\end{eqnarray}
where $p\ge q \ge 1$.   Note that there is a significant difference between system (\ref{lane-emden}) and equation (\ref{4lane}) in the sense that this system has Hamiltonian structure while the equation has gradient structure,  see \cite{dff,dfm,sz} and references therein.  This system has been of great interest at least in the past two decades.   In particular, the Lane-Emden conjecture stating that $u=v=0$  is the only nonnegative solution for this system where $\frac{1}{p+1} +\frac{1}{q+1}>\frac{n-2}{n}$ has been studied extensively and various methods and techniques are developed to tackle this conjecture. Among these methods, Souplet \cite{so} proved the following pointwise inequality for solutions of (\ref{lane-emden}) and then used it to prove  the Lane-Emden conjecture in four dimensions. Note that the particular case $1<p<2$ is done by Phan in \cite{phan}.
\begin{thm} (Souplet \cite{so},  2009) Let $u$ and $v$ be nonnegative solutions of (\ref{lane-emden}).  Then the following inequality holds
\begin{equation}\label{pointlane}
\frac{u^{p+1}}{p+1}\le \frac{v^{q+1}}{q+1}  \ \ \text{in}\ \ \mathbb{R}^n.
\end{equation}
\end{thm}
Applying this theorem, the following pointwise inequality holds for nonnegative solutions of \eqref{4lane}
 \begin{equation}\label{point4lane}
-\Delta u \ge \sqrt\frac{2}{p+1} u^{\frac{p+1}{2}}  \ \ \text{in}\ \ \mathbb{R}^n.
 \end{equation}
   Note also that Phan in \cite{phan}, with similar methods provided in \cite{so}, extended the pointwise inequality (\ref{pointlane}) to nonnegative solutions of the H\'{e}non-Lane-Emden system that is
 \begin{eqnarray}\label{systemhenon}
 \left\{ \begin{array}{lcl}
\hfill -\Delta u&=& |x|^b v^q   \ \ \text{in}\ \ \mathbb{R}^n,\\
\hfill -\Delta v&=& |x|^a u^p   \ \ \text{in}\ \ \mathbb{R}^n,
\end{array}\right.
\end{eqnarray}
where $p\ge q \ge 1$.  Suppose that $0\le a -b\le (n-2)(p-q)$ then
\begin{equation}\label{pointsystemhenon}
|x|^a \frac{u^{p+1}}{p+1}\le |x|^b \frac{v^{q+1}}{q+1}   \ \ \text{in}\ \ \mathbb{R}^n.
 \end{equation}
The standard method to prove a pointwise inequality, as it is used to prove (\ref{pointlane}) and (\ref{pointmod}), is to derive an appropriate equation, call it an auxiliary equation, for the difference function of the right-hand  and the left-hand sides of the inequality.  Then, whenever we have enough decay estimates on solutions of the auxiliary equation,  maximum principles can be applied to prove that the difference function has a fixed sign.     So, the key point here is to manipulate a suitable auxiliary equation.

In a more technical framework,  to construct an auxiliary equation to prove (\ref{pointlane}) and (\ref{point4lane})  a few  positive terms including a gradient term of the form $|\nabla u|^2 u^{t-2}$  for some number $t$ are not considered in \cite{so}.    To be more explicit,   in order to prove (\ref{point4lane}), that is a particular case of (\ref{pointlane}),  the difference function $w(x):=\Delta u+\sqrt\frac{2}{p+1} u^{\frac{p+1}{2}}$ is considered. Straightforward calculations show that the following auxiliary equation holds
\begin{equation}\label{auxw}
\left( \sqrt \frac{2}{p+1}  u^{\frac{1-p}{2}}\right) \Delta w= \Delta u+ \sqrt\frac{2}{p+1} u^{\frac{p+1}{2}}+  \frac{p-1}{2} \frac{|\nabla u|^2}{u}.
\end{equation}
In order to show that $\Delta w$ is nonnegative when $w$ is nonnegative, via maximum principles for the above equation, the gradient term $\frac{|\nabla u|^2}{u}$ is not considered in \cite{so}.  Note that the above equation (\ref{auxw}) implies, in the spirit,  that the gradient term $\frac{|\nabla u|^2}{u}$ should have an impact on the inequality just like the Laplacian operator and the power term $u^{\frac{p+1}{2}}$.    This is our motivation to attempt to include the gradient term in the inequality (\ref{point4lane}) that gives a lower bound on the Laplacian operator.     Let us briefly mention that  Modica in his proof of (\ref{pointmod}) took advantage of similar gradient terms to construct an auxiliary equation.  Following  ideas provided by Modica \cite{m} and Souplet \cite{so}, as we shall see in the proof of Proposition \ref{propwk},  we manage to keep most of the positive terms when looking for an auxiliary equation.

 In this paper, we develop a Moser iteration type argument to
 prove a  lower bound for the negative Laplacian of positive bounded solutions of (\ref{4henon})
 that involves powers of $u$ and  the new term   $\frac{|\nabla u|^2}{u}$ with $\frac{2}{n-4}$ as the coefficient.
 The remarkable point is that the coefficient $\frac{2}{n-4}$ is what we exactly need in the estimate of the scalar curvature for the conformal metric $g = u^{\frac{2}{n-4}}g_0$.
 
      Here is our main result.
   \begin{thm}\label{mainres} Let $u$ be a bounded positive solution of (\ref{4henon}). Then the following pointwise inequality holds
 \begin{equation}\label{newpoint4henon}
-\Delta u \ge \sqrt\frac{2}{(p+1)-c_n} |x|^{\frac{a}{2}} u^{\frac{p+1}{2}} + \frac{2}{n-4} \frac{|\nabla u|^2}{u} \ \ \text{in}\ \ \mathbb{R}^n
 \end{equation}
 where   $c_n:=\frac{8}{n(n-4)}$ and $0\le a\le \inf_{k\ge 0} A_k$ (defined at (\ref{ak})).
  \end{thm}

\begin{remark}  A natural question here is that what  are the best constants in the inequality (\ref{newpoint4henon})?
\end{remark}

 Let us now put the inequality (\ref{newpoint4henon}) in a more geometric text. By the conformal change $g = u^{\frac{4}{n-4}} g_0$ where $g_0$ is the usual Euclidean metric, the new scalar curvature becomes
 $$ S_g= - \frac{4(n-1)}{n-2} u^{-\frac{n+2}{n-4}} \Delta \left(  u^{\frac{n-2}{n-4}} \right). $$
 An immediate consequence of (\ref{newpoint4henon}) is that the conformal scalar curvature is positive. Note that this can not be deduced from the inequality (\ref{point4lane}).

  The idea of proving a lower bound for the negative of Laplacian operator is also used in the context of nonlinear eigenvalue problems to prove certain regularity results, e.g.  see \cite{ceg}.  Similar pointwise inequalities are used to prove Liouville theorems in the notion of stability in \cite{wxy,wy} and references therein as well.    We would like to mention that Gui in \cite{gui} proved a very interesting Hamiltonian identity for elliptic systems  that may be regarded as a generalization of the Modica's inequality.  He used this identity to rigorously analyze    the structure of level curves of saddle solutions of the Allen-Cahn equation as well as Young's Law for the contact angles in triple junction formation.  Note also that as it is shown by Farina in \cite{far} for the Ginzburg-Landau system, the analog of Modica's estimate  is false for systems in general.   We refer  interested readers to \cite{ali} for a review of this topic  and  to \cite{fg2} for  De Giorgi type results for systems.

  Here is the organization of the paper. In Section \ref{secEst}, we provide certain standard elliptic estimates that are consequences of Sobolev embeddings and the regularity theory. Then, in Section \ref{secIter} we develop a Moser iteration type argument, following ideas provided by Modica \cite{m} and Souplet in \cite{so}.  Finally, in Section \ref{secapp}, we first  give a certain maximum principle type argument for a quasilinear equation that arises in the Moser iteration process. Then we apply the estimates and methods developed in former sections.    We suggest to ignore the weight function $|x|^a$ in (\ref{4henon}) when  reading the paper for the first time.


 \section{Technical elliptic estimates}\label{secEst}

 In this section, we provide some elliptic decay estimates that we use frequently later in the proofs.  Deriving the right decay estimates for solutions of (\ref{4henon}) play a fundamental role in the most our proofs.    Similar estimates have been also used in the literature to construct Liouville theorems and regularity results. We refer the interested readers  to \cite{f,fg,phan,so,ps}. We start with the following standard estimate. 

 \begin{lemma}\label{1bound} ($L^p$-estimate on $B_R$) Suppose that $u$ is a nonnegative solution of (\ref{4henon}) then for any $R>1$ we have
   \begin{equation*}
\int_{B_R} |x|^a u^p  \le C \ R^{n-\frac{4p+a}{p-1}},
 \end{equation*}
 where $C=C(n,p,a)>0$ is independent from $R$. 
\end{lemma}
\noindent\textbf{Proof:} Consider the following test function $\phi_R\in C^4_c(\mathbb{R}^n)$ with $0\le\phi_R\le1$;
 $$\phi_R(x)=\left\{
                      \begin{array}{ll}
                        1, & \hbox{if $|x|<R$;} \\
                        0, & \hbox{if $|x|>2R$;} 
                                                                       \end{array}
                    \right.$$
where $|| D^{i} \phi_R||_{\infty} \le \frac{C}{R^{i}}$ where $1\le i \le 4$. For fixed $m\ge 2$, we have  
$$|\Delta^2 \phi^m_R(x)|\le  \left\{
                      \begin{array}{ll}
                        0, & \hbox{if $|x|<R$ or $|x|>2R$;} \\
                         C R^{-4} \phi^{m-4}_R , & \hbox{if $R<|x|<2R$;} 
                                                                       \end{array}
                    \right.$$
 where $C>0$ is independent from $R$.  For $m\ge 2$, multiply the equation by $\phi^m_R$ and integrate to get 
\begin{eqnarray*}
 \int_{B_{2R}}     |x|^{a}     u^p  \phi^m_R &= & \int_{B_{2R}}  \Delta^2 u   \phi^m_R\\
&=& \int_{B_{2R}} u \Delta^2 \phi^m_R   \le C  R^{-4}\int_{B_{2R}\setminus B_{R}} u\phi^{m-4}_R.
\end{eqnarray*}
Applying H\"{o}lder's inequality we get 
\begin{eqnarray*}
\int_{B_{2R}}      |x|^{a}     u^p  \phi^m_R & \le  & C \ R^{-4}  \left(         \int_{B_{2R}\setminus B_{R}}  |x|^{\frac{-a}{p }p'    }          \right)^{\frac{1}{p'}}   \left(     \int_{B_{2R}\setminus B_{R}}      |x|^{a}    u^p \phi^{(m-4)p}_R   \right)^{1/p}\\
& \le  & C\  R^{  (n-\frac{a}{p}p')\frac{1}{p'}  -4 }  \left(     \int_{B_{2R}\setminus B_{R}}       |x|^{a}     u^p \phi^{(m-4)p}_R   \right)^{1/p},
\end{eqnarray*}
where $p'=\frac{p}{p-1}$.  Set $m=(m-4)p$ that gives $m=\frac{4p}{p-1}$ to get 
\begin{equation*}
\int_{B_{2R}}      |x|^{a}     u^p  \phi^m_R \le    C\  R^{  (n-\frac{a}{p}p')\frac{1}{p'}  -4 }  \left(     \int_{B_{2R}}       |x|^{a}     u^p \phi^{m}_R   \right)^{1/p}.
\end{equation*}
Therefore, 
\begin{equation*}
\int_{B_{2R}}      |x|^{a}     u^p  \phi^m_R \le  C\  R^{  (n-\frac{a}{p}p')  -4p'}  .
\end{equation*}
This finishes the proof. 

\hfill $\Box$

From the H\"{o}lder's inequality we get the following.
    \begin{cor}\label{u} Under the same assumptions as Lemma \ref{1bound}. The following estimate holds
   $$ \int_{B_R\setminus B_{R/2}}  u \le C R^{n -\frac{a+4}{p-1}}$$
  where $C=C(n,p,a)>0$  is independent from $R$.
    \end{cor}

We now show that the operator $-\Delta u$ has a sign. Then, we apply this to provide various elliptic estimates for derivatives of $u$. In addition, later on this helps us to start an iteration argument.   

 \begin{prop}\label{prop0}
  Let $u$ be a positive solution of (\ref{4henon}). Then, $-\Delta u\ge 0$ in $\mathbb {R}^n$.
 \end{prop}

\noindent\textbf{Proof:}  Let $v=-\Delta u$. Ideas and methods applied in this proof are strongly motivated by the ones given in \cite{wx}.    Suppose that there is $x_0\in\mathbb R^n$ such that $v(x_0)<0$. Without loss of generality we take $x_0=0$, i. e.  in case of $x_0\neq 0$ set $\omega(x)=v(x+x_0)$ and apply the same argument.  We use the notation $\bar f(r)=\frac{1}{|\partial B_r|} \int_{\partial B_r} f dS$ as the average of function $f(x)$ on the boundary of $B_r$.  We  refer interested readers to \cite{n} regarding the average function.   Applying the H\"{o}lder's inequality
\begin{eqnarray}\label{bar}
 \left\{ \begin{array}{lcl}
\hfill -\Delta_r \bar u (r)&=& \bar v (r)  \ \ \text{in}\ \ \mathbb{R},\\   
\hfill -\Delta_r \bar v(r) &\ge & r^a (\bar u)^p   \ \ \text{in}\ \ \mathbb{R},
\end{array}\right.
\end{eqnarray}
where $\Delta_r$  is the Laplacian operator in the polar coordinates, i.e. $$\Delta_r \bar f(r)= r^{1-n} ( r^{n-1} \bar f'(r) )'.$$ It is straightforward to see that 
 $$\bar v'(r)=\frac{1}{|\partial B_r|} \int_{B_r} \Delta v=-\frac{1}{|\partial B_r|} \int_{B_r} |x|^a u^p \le 0.$$ 
Therefore, $\bar v(r)\le \bar v(0)<0$ for $r>0$. Similarly for $\bar u'(r)$ we have
\begin{eqnarray*}
\bar u'(r) &=&- \frac{1}{|\partial B_r|} \int_{B_r} v = - r^{1-n} \int_0^r s^{n-1} \bar v(s) ds
\\&\ge&  - \bar v(0)  r^{1-n} \int_0^r  s^{n-1} ds = - \frac{\bar v(0)}{n} r.
\end{eqnarray*}
From this for any $r\ge r_0$ we get 
\begin{eqnarray}\label{bar2}
\bar u(r)  \ge \alpha r^2 ,
\end{eqnarray}
where $\alpha=- \frac{\bar v(0)}{2n} >0$.  We now have a lower bound on $\bar u(r) $.  
Instead suppose that the following more general lower bound holds on $\bar u(r) $,
\begin{eqnarray}\label{bar2}
\bar u(r)  \ge \frac{\alpha^{p^k}}{\beta^{s_k}} r^{t_k}  \ \ \text{for} \ \ r\ge r_k,
\end{eqnarray}
where  $s_0:=0$, $t_0:=2$, $\alpha:=- \frac{\bar v(0)}{2n}>0$  and $\beta:=2p+a+n+4>0$.      Note that system (\ref{bar}) makes a relation between two functions $\bar u(r) $ and $\bar v(r) $. Therefore, the lower bound on   $\bar u(r) $ forces an upper bound on $\bar v(r) $ and vice versa.   In the light of this fact,  we can construct an iteration argument to improve the bound (\ref{bar2}).  Integrating the second equation of (\ref{bar}) over $[r_k,r]$ when $r\ge r_k$ we get 
\begin{eqnarray*}
r^{n-1} \bar v'(r)& \le& r_k^{n-1} \bar v'(r_k)-  \frac{\alpha^{p^{k+1}}}{\beta^{p s_k}} \int_{r_k}^{r} s^{n-1+a+pt_k} ds\\&\le&-  \frac{\alpha^{p^{k+1}}}{\beta^{p s_k} (pt_k+n+a)} (r^{pt_k+n+a}-r_k^{pt_k+n+a}) \ \ \text{since} \ \ \bar v' <0. 
\end{eqnarray*}
Therefore $\bar v'(r) \le -  \frac{\alpha^{p^{k+1}}}{\beta^{p s_k} (pt_k+n+a)} (r^{pt_k+a+1}-r_k^{pt_k+a+1}) \ \ \text{for all}  \ \   r \ge r_k$  that is 
\begin{equation*}
 \bar v'(r) \le -  \frac{\alpha^{p^{k+1}}}{2\beta^{p s_k} (pt_k+n+a)} r^{pt_k+a+1} \ \ \text{for all}  \ \   r \ge    2^{\frac{1}{pt_k+a+1}}r_k. 
\end{equation*}
Integrating the last inequality over $[2^{\frac{1}{pt_k+a+1}}r_k,r]$ when $r\ge 2^{\frac{1}{pt_k+a+1}}r_k=\tilde r_k$,  we obtain
\begin{equation*}
 \bar v(r) \le  \bar v(\tilde r_k) -  \frac{\alpha^{p^{k+1}}}{2\beta^{p s_k} T_{k,n,a,p}} (r^{pt_k+a+2}-\tilde r_k^{pt_k+a+2} ),
\end{equation*}
where $T_{k,n,a,p}:=(pt_k+n+a)(pt_k+2+a)$. By similar discussions and by taking $r $ large enough,  that is $r\ge    2^{\frac{1}{pt_k+a+1}}  2^{\frac{1}{pt_k+a+2}} r_k ={\tilde {\tilde r}}_k$,  we end up with 
\begin{equation}\label{vr}
 \bar v(r) \le -  \frac{\alpha^{p^{k+1}}}{4\beta^{p s_k} T_{k,n,a,p}} r^{pt_k+a+2}.
\end{equation}
Applying (\ref{vr}) and integrating equation (\ref{bar}) again over $[{\tilde {\tilde r}}_k,r]$ when  $r\ge {\tilde {\tilde r}}_k$,  we have 
\begin{eqnarray*}
r^{n-1} \bar u'(r)&=& \tilde r_k^{n-1} \bar u'(\tilde r_k)  - \int_{\tilde r_k}^{r} s^{n-1} \bar v(s) ds
\\&\ge&  \frac{\alpha^{p^{k+1}}}{4\beta^{p s_k} T_{k,n,a,p}} \int_{\tilde r_k}^r   s^{pt_k+a+n+1} ds.
\end{eqnarray*}
Therefore, the following new lower bound on $\bar u(r)$ holds
\begin{eqnarray*}
\bar u(r)&\ge& \frac{\alpha^{p^{k+1}}}{2^4 \beta^{p s_k} \tilde T_{k,n,a,p}}   r^{pt_k+a+n+4},
\end{eqnarray*}
where $$r\ge    2^{\frac{1}{pt_k+a+3}}  2^{\frac{1}{pt_k+a+4}} {\tilde {\tilde r}}_k=2^{ \sum_{i=1}^4 \frac{1}{pt_k+a+i} }r_k , $$ 
and 
\begin{eqnarray*}
\tilde T_{k,n,a,p}&=&  (pt_k+n+a+2) (pt_k+4+a) T_{k,n,a,p}\\&=&(pt_k+n+a)(pt_k+2+a) (pt_k+n+a+2) (pt_k+4+a)\\&\le& (pt_k+n+a+4)^4 .
\end{eqnarray*}
 We now modify this estimate to make the coefficients similar to (\ref{bar2}).  After simplifying  we get
\begin{equation}\label{bar3}
\bar u(r)  \ge  \frac{\alpha^{p^{k+1}}}{\beta^{p s_k} M_k}   r^{pt_k+a+4} \ \ \ \ \text{for} \ \  r\ge 2^{ \frac{4}{pt_k+a+1}} r_k ,
\end{equation}
where $M_k:=2^4 (pt_k+n+a+4)^4$.  In what follows,  we put an upper bound on $M_k$ that is expressed as a power of  $\beta$. Note that 
\begin{eqnarray*}
\frac{1}{2}\sqrt[4]{M_{k+1}}&=&pt_{k+1}+n+a+4= p(pt_{k}+n+a+4)++n+a+4 \\&\le& (pt_{k}+n+a+4)(p+1)=\frac{p+1}{2}\sqrt[4]{M_{k}}. 
\end{eqnarray*}
From this we have $M_{k+1} \le (p+1)^4M_k$ and therefore  $M_k \le (p+1)^{4k}M_0$ where    $M_0=2^4(2p+n+a+4)^4$ because $t_0=2$. Since the constant $\beta$ is defined as $\beta=2p+n+a+4$, we get the following bound
\begin{equation}\label{mk}
M_k \le \beta^{4k+4}.
\end{equation}
 From this, (\ref{bar2}) and (\ref{bar3}) and to complete the iteration process, we set 
\begin{eqnarray}\label{tk}
t_{k+1}&:=&pt_k+a+4 \ \ \text{for } \ t_0=2,\\ \label{sk}
s_{k+1}&:=& p s_k+ 4k+4 \ \ \text{for} \ s_0=0,
\end{eqnarray}
and therefore, 
\begin{eqnarray}\label{bar4}
\bar u(r)  \ge \frac{\alpha^{p^{k+1}}}{\beta^{s_{k+1}}} r^{t_{k+1}}
  \ \ \text{for} \ \ r\ge r_{k+1},
\end{eqnarray}
where $ r_{k+1}:=2^{ \frac{4}{pt_k+a+1}} r_k \ge 2^{ \sum_{i=1}^4 \frac{1}{pt_k+a+i} }r_k .$  
By direct calculations on these recursive sequences we get the explicit sequences 
\begin{eqnarray*}
t_k&=& \frac{2p^{k+1} +(a+2) p^k-(a+4)}{p-1}, \\ 
s_k&=&\frac{4p^{k+1}-4p(k+1)+4k}{(p-1)^2} , \\
r_{k}&=&2^{ \sum_{i=0}^{k-1} \frac{4}{pt_i+a+1}} r_0  \le 2^{ \sum_{i=0}^{\infty} \frac{4}{pt_i+a+1}} r_0 =:r^* <\infty.
\end{eqnarray*}
Set $R:=\beta^{\frac{2}{p-1}} M$ where  $M=\max\{\alpha^{-1},m\}$ when $m>1$ is large enough to make sure $m \beta^{\frac{2}{p-1}} \ge r^* $.  Therefore,  $R \ge r^* \ge r_k$ for any $k$ and we have 
$$ 
\bar u(R) \ge  M^{t_k-p^k}  \beta^{ \frac{2t_k}{p-1} -s_k }.
$$
 If we take $k$ large enough, e.g.  $k \ge \frac{\ln (a+4)-\ln(a+2)}{\ln p}$,  then $t_k>p^k$. The fact that  $M>1$, gives us
\begin{equation*}
\bar u(R) \ge  \beta^{ \frac{2t_k}{p-1} -s_k } = \beta^{   \frac{   2(a+2) p^k +4k(p-1) +4p-2(a+4)  }{(p-1)^2}  } .
\end{equation*}
Since we have assumed that $a+2>0$ and $\beta>1$,  we get $\bar u(R)\to\infty$  as $k\to\infty$. Note that $0<R<\infty$ is independent from $k$.  This finishes the proof.

    \hfill $\Box$

We now apply Proposition \ref{prop0} to conclude that $-\Delta u \ge 0$ and therefore we can consider equation (\ref{4henon}) as a special case of the H\'{e}non-Lane-Emden equation. 

 \begin{lemma}\label{1boundLap} ($L^1$-estimates on $B_R$) Suppose that $u$ is a nonnegative solution of (\ref{4henon}) then for any $R>1$ we have
   \begin{equation*}
\int_{B_R} |\Delta u| \le C R^{n-\frac{2p+2+a}{p-1}} ,\\
 \end{equation*}
 where $C=C(n,p,a)>0$  is independent from $R$.
 \end{lemma}
 \noindent\textbf{Proof:} Set $v=-\Delta u$. From Proposition \ref{prop0} we know that $v\ge0$. Therefore the pair  $(u,v)$ satisfies the following system 
 \begin{eqnarray}
\label{mainbound}
 \left\{ \begin{array}{lcl}
\hfill -\Delta u&=& v   \ \ \text{in}\ \ \mathbb{R}^n,\\   
\hfill -\Delta v&=& |x|^{a}u^p   \ \ \text{in}\ \ \mathbb{R}^n,
\end{array}\right.
  \end{eqnarray}
  that is a particular case of the H\'{e}non-Lane-Emden system. From the estimates provided in \cite{fg} as Lemma 2.1 we get the desired result.

\hfill $\Box$

\begin{lemma}\label{interp}
(An interpolation inequality on $B_R$)  Let  $R>1$ and $z\in W^{2,1}(B_{2R})$.  Then
$$\int_{B_R\setminus B_{R/2}}     |  D z |\le C R \int_{B_{2R}\setminus B_{R/4}}   |\Delta z|  +    C R^{-1}   \int_{B_{2R}\setminus B_{R/4}}    |z|   ,$$
where $C=C(n)>0$ is independent from $R$.
\end{lemma}

\begin{cor}\label{corgrad} Under the same assumptions as Lemma \ref{1bound}. The following estimate holds.
$$\int_{B_R\setminus B_{R/2}}     |  D u |\le C R^{n-\frac{p+3+a}{p-1}} ,$$
where $C=C(n,p,a)>0$ is independent from $R$.
\end{cor}

\begin{lemma} \label{ellip}
($L^\tau$-estimate on $B_R$)  Let $1< \tau<\infty$ and $z\in W^{2,\tau}(B_{2R})$. Then,
$$\int_{B_R\setminus B_{R/2}}   |  D^{2} z |^\tau\le C \int_{B_{2R}\setminus B_{R/4}}   |\Delta z|^\tau  +  C R^{-2\tau}    \int_{B_{2R}\setminus B_{R/4}}     |z|^\tau ,$$
where $C=C(n,\tau)>0$  does not depend on $R$.
\end{lemma}

 \begin{lemma}\label{2bound} ($L^2$-estimates on $B_R$) Suppose that $u$ is a bounded nonnegative solution of (\ref{4henon}) then for any $R>1$ we have
  \begin{equation}\label{stepfin}
\int_{B_R} |\Delta u|^2 \le C  \int_{B_{2R}} |x|^a u^{p+1} + C R^{-2}  \int_{B_{2R}}   |\Delta u| + C R^{-4}  \int_{B_{2R}\setminus B_{R}}   u,
\end{equation}
 where $C=C(n,p,a)>0$   does not depend on $R$.
\end{lemma}
\noindent\textbf{Proof:} We proceed in two steps.

\noindent Step 1.  Multiply the both sides of equation (\ref{4henon}) with $u \phi^2$ where $\phi\in C_c^\infty(\mathbb{R}^n)\cap[0,1]$ is a test function. Then, doing the integration by parts, we get
\begin{eqnarray*}
\int_{\mathbb{R}^n} |\Delta u|^2 \phi^2 &=& \int_{\mathbb{R}^n} |x|^a u^{p+1} \phi^2 - 4 \int_{\mathbb{R}^n} \Delta u \nabla u\cdot \nabla \phi \phi - \int_{\mathbb{R}^n} u \Delta u \left(2 |\nabla\phi|^2+2 \phi\Delta \phi\right) 
 \\&\le&     \int_{\mathbb{R}^n} |x|^a u^{p+1} \phi^2 + \delta \int_{\mathbb{R}^n} |\Delta u|^2 \phi^2 + C(\delta)  \int_{\mathbb{R}^n} |\nabla u|^2 |\nabla \phi|^2 \\&&+ C  \int_{\mathbb{R}^n} |\Delta u | \left( |\nabla\phi|^2+|\Delta \phi|\right),
\end{eqnarray*}
for some constant $C>0$. Here, we have used the Cauchy's inequality for $0<\delta <1$. Therefore, if we set $\phi$ to be the standard test function that is $\phi=1$ in $B_R$ and $\phi=0$ in $\mathbb{R}^n\setminus B_{2R}$ when $||D_x^i \phi ||_{L^{\infty}(B_{2R}\setminus B_R)} \le C R^{-i}$ for $i=1,2$, then we get
\begin{equation}\label{step1}
\int_{B_R} |\Delta u|^2 \le  \int_{B_{2R}} |x|^a u^{p+1} + C R^{-2}  \int_{B_{2R}\setminus B_{R}}   |\nabla u|^2 + C R^{-2}  \int_{B_{2R}\setminus B_{R}}   |\Delta u|,
\end{equation}
where $C=C(n,p,a)>0$   does not depend on $R$.

 \noindent Step 2.  Multiply the both sides of  $-\Delta u=v$ with $u \phi^2$ where $\phi$ is the same test function as Step 1.  Again doing integration by parts we get
 \begin{eqnarray*}
\int_{\mathbb{R}^n} |\nabla u|^2 \phi^2 &=& \int_{\mathbb{R}^n} uv \phi^2 - 2 \int_{\mathbb{R}^n} u\nabla u\cdot \nabla \phi \phi  
\\&\le&    \int_{\mathbb{R}^n} uv \phi^2  + \delta \int_{\mathbb{R}^n} |\nabla u|^2 \phi^2 + C(\delta)  \int_{\mathbb{R}^n}  |\nabla \phi|^2 u^2 ,
\end{eqnarray*}
where we have also used the Cauchy's inequality for $0<\delta <1$.   So,
\begin{equation}\label{step2}
\int_{B_R} |\nabla u|^2 \le C   \int_{B_{2R}}   |\Delta u| + C R^{-2}  \int_{B_{2R}\setminus B_{R}}   u,
\end{equation}
where we have used the boundedness of $u$. From (\ref{step1}) and (\ref{step2}) we get
\begin{equation}\label{stepfin}
\int_{B_R} |\Delta u|^2 \le   \int_{B_{2R}} |x|^a u^{p+1} + C R^{-2}  \int_{B_{2R}}   |\Delta u| + C R^{-4}  \int_{B_{2R}\setminus B_{R}}   u.
\end{equation}
This completes the proof. 

\hfill $\Box$

We now apply Lemma \ref{1bound}, Lemma \ref{2bound} and Corollary \ref{u} to get the following.
    \begin{cor}\label{Delta2} Suppose that the assumptions of Lemma \ref{1bound} hold. Moreover, let $u$ be bounded then
       \begin{equation}\label{Delta2R}
        \int_{B_R} |\Delta u|^2 \le  C R^{n-\frac{4p+a}{p-1}},
        \end{equation}
        where $C=C(n,p,a)>0$ is independent from $R$.
    \end{cor}

\begin{lemma} (Sobolev inequalities on the sphere $S^{n-1}$) \label{sobolev}
Let $n\ge2$,  integer $i\ge 1$ and $1<t<\tau\le\infty$. For $z\in W^{i,t}(S^{n-1})$, the following estimate holds
     $$||z||_{L^\tau(S^{n-1})}\le C || D_\theta^i  z||_{L^t(S^{n-1})}  + C || z  ||_{L^1(S^{n-1})} ,$$ where

       $$\left\{
                      \begin{array}{ll}
                        \frac{1}{\tau}= \frac{1}{t}-\frac{i}{n-1}, & \hbox{if $it+1<n$,} \\
                        \tau=\infty, & \hbox{if $it+1>n$,}
                                                                       \end{array}
                    \right.$$
                     and $C=C(i,t,n,\tau)>0$.
\end{lemma}


 \section{Developing the iteration argument}\label{secIter}

 In this section, we develop a counterpart of the Moser iteration argument \cite{mos} for solutions of (\ref{4henon}).  We define a sequence of functions $(w_k)_{k=-1}$ of the form $$  w_k:=\Delta u+\alpha_k |\nabla u|^2 (u+\epsilon)^{-1} + \beta_k |x|^{\frac{a}{2}} u^{\frac{p+1}{2}}
$$ where $\alpha_k$ and $\beta_k$ are certain nondecreasing sequences of nonnegative numbers where $\alpha_{-1}=\beta_{-1}=0$.

 Assuming that $w_k\le 0$,  that is essentially a lower bound on the negative Laplacian operator,  holds  we construct a differential inequality for $w_{k+1}$ where $\alpha_{k+1}\ge \alpha_k$ and $\beta_{k+1}\ge \beta_k$.  Then, applying certain maximum principle type arguments, we show that $w_{k+1}\le 0$. Note that  $w_{k+1}\le 0$ is stronger than $w_k\le 0$,  because it forces a stronger lower bound on the negative of Laplacian operator.

We start with proving that $w_{-1}$, which is the Laplacian operator of $u$, is nonpositive, see Proposition \ref{prop0}.   Then,  using this fact and applying (\ref{systemhenon}) and (\ref{pointsystemhenon}) when $q=1$ and $b=0$,  we get the following inequality for nonnegative solutions of the fourth order H\'{e}non equation (\ref{4henon})
\begin{equation}\label{point4henon}
-\Delta u \ge \sqrt\frac{2}{p+1} |x|^{\frac{a}{2}} u^{\frac{p+1}{2}}  \ \ \text{in}\ \ \mathbb{R}^n,
 \end{equation}
 where $0\le a\le (n-2)(p-1)$. Inequality (\ref{point4henon}) is the first step of the iteration argument  meaning that $w_0 \le 0$ for $\alpha_0=0$ and $\beta_0=\sqrt\frac{2}{p+1}$.   

 We now perform the iteration argument.

\begin{prop}\label{propwk}
Let $u$ be a positive classical solution of (\ref{4henon}). Suppose that $(\alpha_k)_{k=0}$ and $(\beta_k)_{k=0}$ are sequences of numbers.  Define the following sequence of functions
 \begin{equation}\label{wk}
  w_k:=\Delta u+\alpha_k |\nabla u|^2 (u+\epsilon)^{-1} + \beta_k |x|^{\frac{a}{2}} u^{\frac{p+1}{2}},
   \end{equation}
where $\epsilon=\epsilon(k)$ is a positive constant.   Suppose that $w_k\le 0$, then $w_{k+1}$ satisfies the following differential inequality
 \begin{eqnarray}\label{wk+1}
 &&\Delta w_{k+1}-2 \alpha_{k+1} (u+\epsilon)^{-1} \nabla u\cdot\nabla w_{k+1}\\ &&\nonumber +\alpha_{k+1} w_{k+1} (u+\epsilon)^{-2} |\nabla u|^2 -  \frac{\beta_{k+1}(p+1)}{2} u^{\frac{p-1}{2}} |x|^{\frac{a}{2}} w_{k+1}\\ &&\nonumber \ge
  I^{(1)}_{\epsilon,\beta_{k} }   |x|^a u^p + \alpha_{k+1}  I^{(2)}_{\alpha_{k}} |\nabla u|^4 (u+\epsilon)^{-3} + I^{(4)}_{a,\alpha_{k},\beta_{k}} |x|^{a-2} u^{\frac{p+1}{2}}
  \\&&\nonumber+ I^{(3)}_{\epsilon,\alpha_{k},\beta_{k}}  |x|^a u^{\frac{p+1}{2}}  \left| \frac{\nabla u}{u} +\frac{  a\beta_{k+1}   \left(\frac{p+1}{2}- \alpha_{k+1} \frac{u}{u+\epsilon}\right) }{ 2 I^{(3)}_{\epsilon,\alpha_{k},\beta_{k}} }\frac{x}{|x|^2} \right|^2,
  \end{eqnarray}
where
 \begin{eqnarray*}
 I^{(1)}_{\epsilon, \alpha_k,\beta_{k}}&: =& 1- \frac{p+1}{2} \beta_{k+1}^2 +\frac{2}{n} \alpha_{k+1} \beta_k^2 \frac{u}{u+\epsilon},
 \\
  I^{(2)}_{\alpha_{k}} &:= &\frac{2}{n} (\alpha_{k+1} +\alpha_k+1)^2 -2\alpha_{k+1} (\alpha_{k+1}+1)+\alpha_{k+1},
  \\
   I^{(3)}_{\epsilon,\alpha_{k},\beta_{k}}&: = &  \frac{4}{n} \alpha_{k+1} \beta_k (\alpha_{k+1} +\alpha_k+1) \frac{u^2}{(u+\epsilon)^2} +\beta_{k+1} \alpha_{k+1} \frac{u^2}{(u+\epsilon)^2}
   \\&& -(p+1) \beta_{k+1} \alpha_{k+1} \frac{u}{(u+\epsilon)} +\frac{p+1}{2}\left(\frac{p-1}{2}-\alpha_{k+1} \frac{u}{(u+\epsilon)} \right ) \beta_{k+1},
   \\
    I^{(4)}_{a,\epsilon, \alpha_{k},\beta_{k}}&: = & \frac{a}{2}\beta_{k+1} (n+\frac{a}{2} -2) - \frac{a^2 \beta^2_{k+1}   \left(\frac{p+1}{2}- \alpha_{k+1} \frac{u}{u+\epsilon}\right) ^2   }{   4  I^{(3)}_{\epsilon,\alpha_{k},\beta_{k}}  } .
 \end{eqnarray*}
\end{prop}

\noindent\textbf{Proof:} For the sake of simplicity in calculations, set $b:=\frac{a}{2}$ and $q:=\frac{p+1}{2}$.   From (\ref{wk}) the function $w_{k+1}$ is defined as $$ w_{k+1}:=\Delta u+\alpha_{k+1} |\nabla u|^2 (u+\epsilon)^{-1} + \beta_{k+1} |x|^{b} u^q.$$  
Taking Laplacian of $w_{k+1}$  and using equation (\ref{4henon}) we get
\begin{eqnarray}\label{wk++1}
\Delta w_{k+1}&=& \Delta^2 u+ \alpha_{k+1}  \Delta (|\nabla u|^2  (u+\epsilon)^{-1}  ) +\beta_{k+1} \Delta (|x|^b u^q)
\\ \nonumber &=&|x|^a u^p+I+J,
\end{eqnarray}
where $I:= \alpha_{k+1}  \Delta (|\nabla u|^2(u+\epsilon)^{-1}) $ and $J:=\beta_{k+1} \Delta (|x|^b u^q)$. In what follows, we simplify $I$ and $J$ as well as finding lower bounds for these terms.     We start with $J$ that is
 \begin{eqnarray*}
 \frac{J}{\beta_{k+1}}&=& \Delta (|x|^b u^q)= \Delta |x|^b u^q +\Delta u^q |x|^b +2\nabla |x|^b\cdot \nabla u^q \\ &=& b (n+b-2) |x|^{b-2} u^q + q(q-1) |x|^b u^{q-2}  |\nabla u|^2\\& &+q |x|^b u^{q-1} \Delta u  +2b q |x|^{b-2} u^{q-1} \nabla u\cdot x.
 \end{eqnarray*}
 From the definition of $w_{k+1}$, we have
 \begin{equation}\label{uwk+1}
  \Delta u=w_{k+1}-\alpha_{k+1} |\nabla u|^2 (u+\epsilon)^{-1} -\beta_{k+1} |x|^{b} u^q.
  \end{equation}
   Substitute this into the last equation to simplify $J$ as
 \begin{eqnarray}\label{J}
 \frac{J}{\beta_{k+1}}&=&  q u^{q-1} |x|^b w_{k+1}-q\beta_{k+1} u^{2q-1} |x|^{2b}  \\&& \nonumber  +
 \left(q(q-1)-q\alpha_{k+1}\frac{u}{u+\epsilon}\right) |x|^b u^{q-2}  |\nabla u|^2 
 \\&& \nonumber +b (n+b-2) |x|^{b-2} u^q  +2b q |x|^{b-2} u^{q-1} \nabla u\cdot x.
 \end{eqnarray}
We now simplify $I$ as what follows,
\begin{eqnarray*}
\frac{I}{\alpha_{k+1}}&=& \Delta (|\nabla u|^2 (u+\epsilon)^{-1} )= \sum_{i,j}  \partial_{jj}(u^2_i (u+\epsilon)^{-1} )\\&=& 2 (u+\epsilon)^{-1}  \sum_{i,j} (\partial_{ij}u)^2  +2 (u+\epsilon)^{-1}  \nabla u\cdot \nabla \Delta u -4 (u+\epsilon)^{-2}  \sum_{i,j} \partial_i u \partial_j u \partial_{ij} u \\&&- |\nabla u|^2 (u+\epsilon)^{-2}  \Delta u + 2 |\nabla u|^4 (u+\epsilon)^{-3} .
 \end{eqnarray*}
Again substituting (\ref{uwk+1}) into the term $2(u+\epsilon)^{-1}   \nabla u\cdot \nabla \Delta u$ appeared above, we get
\begin{eqnarray*}
\frac{I}{\alpha_{k+1}}&=&2 (u+\epsilon)^{-1}  \sum_{i,j} (\partial_{ij}u)^2   -4 (u+\epsilon)^{-2}  \sum_{i,j} \partial_i u \partial_j u \partial_{ij} u \\&&+ 2 |\nabla u|^4 (u+\epsilon)^{-3} - |\nabla u|^2 (u+\epsilon)^{-3} \Delta u \\&&+2 (u+\epsilon)^{-1}  \nabla u\cdot\nabla w_{k+1}-2\alpha_{k+1} (u+\epsilon)^{-1}  \nabla u\cdot\left(|\nabla u|^2(u+\epsilon)^{-1} \right) \\&& -2 \beta_{k+1} (u+\epsilon)^{-1}  \nabla u\cdot \nabla\left( |x|^b u^q \right).
 \end{eqnarray*}
Then, collecting the similar terms we obtain
\begin{eqnarray*}
&& \frac{I}{\alpha_{k+1}} -2 (u+\epsilon)^{-1} \nabla u\cdot\nabla w_{k+1} 
= 2(u+\epsilon)^{-1} \sum_{i,j} (\partial_{ij}u)^2  
\\&& -4(\alpha_{k+1}+1) (u+\epsilon)^{-2} \sum_{i,j} \partial_i u \partial_j u \partial_{ij} u 
\\&&+ 2(\alpha_{k+1}+1) |\nabla u|^4 (u+\epsilon)^{-3}- |\nabla u|^2 (u+\epsilon)^{-2} \Delta u\\&& -2 \beta_{k+1}b |x|^{b-2} (u+\epsilon)^{-1} u^{q} \nabla u\cdot x 
\\&&-2\beta_{k+1} q |x|^{b} u^{q-1} (u+\epsilon)^{-1} |\nabla u|^2.
 \end{eqnarray*}
Completing the square we get
\begin{eqnarray}\label{Isquare}
&& \frac{I}{\alpha_{k+1}}-2 (u+\epsilon)^{-1} \nabla u\cdot\nabla w_{k+1} 
\\&&\nonumber= 2(u+\epsilon)^{-1} \sum_{i,j} \left(  \partial_{ij}u -(\alpha_{k+1}+1) (u+\epsilon)^{-1}   \partial_i u \partial_j u \right)^2
\\&& \nonumber  -2 \alpha_{k+1} (\alpha_{k+1}+1) |\nabla u|^4 (u+\epsilon)^{-3} - |\nabla u|^2 (u+\epsilon)^{-2} \Delta u
 \\&& \nonumber  -2 \beta_{k+1}b |x|^{b-2} (u+\epsilon)^{-1} u^{q} \nabla u\cdot x -2\beta_{k+1} q |x|^{b} u^{q-1} (u+\epsilon)^{-1} |\nabla u|^2.
 \end{eqnarray}
Note that for any $n\times n$ matrix $A=(a_{i,j})$ the Hilbert-Schmidt norm is defined by $|| A||_2=\sqrt {\sum_{i,j}  | a_{i,j} |^2} =\sqrt {\text{trace} (A A^*)}$, where  $A^*$ denotes the conjugate transpose of $A$. From the Cauchy-Schwarz' inequality,  the following inequality holds,
\begin{equation}\label{matrix}
|\text{trace} \ A|^2=|(A,I)|^2 \le ||A||_2^2 ||I||_2^2=n \sum_{i,j}  | a_{i,j} |^2 .
\end{equation}
 Set $a_{i,j}:=   \partial_{ij}u -(\alpha_{k+1}+1) (u+\epsilon)^{-1}   \partial_i u \partial_j u $ in (\ref{matrix}) to get
 $$\sum_{i,j}^{n} \left(   \partial_{ij}u -(\alpha_{k+1}+1) (u+\epsilon)^{-1}  \partial_i u \partial_j u \right)^2 \ge \frac{1}{n} \left(  \Delta u -(\alpha_{k+1}+1) (u+\epsilon)^{-1}   |\nabla u|^2  \right)^2 .
 $$
 From this lower bound for the Hessian and (\ref{Isquare}), we get
\begin{eqnarray}\label{Ialpha}
&&\frac{I}{\alpha_{k+1}} -2 (u+\epsilon)^{-1}  \nabla u\cdot\nabla w_{k+1}  \\ & &\ge  \nonumber  \frac{2}{n}(u+\epsilon)^{-1}  \left(  \Delta u -(\alpha_{k+1}+1)(u+\epsilon)^{-1}   |\nabla u|^2  \right)^2 \\ && \nonumber  -2 \alpha_{k+1} (\alpha_{k+1}+1) |\nabla u|^4 (u+\epsilon)^{-3}  - |\nabla u|^2 (u+\epsilon)^{-2} \Delta u +T_{k},
 \end{eqnarray}
 where $$T_k:=-2 \beta_{k+1}b |x|^{b-2} (u+\epsilon)^{-1} u^{q} \nabla u\cdot x -2\beta_{k+1} q |x|^{b} u^{q-1} (u+\epsilon)^{-1} |\nabla u|^2.$$
  Note also that from the assumption $w_k \le 0$ we have this upper bound on the Laplacian operator,  $\Delta u\le -\alpha_k |\nabla u|^2 (u+\epsilon)^{-1} - \beta_k |x|^{b} u^q$.  Elementary calculations show that if $t\le t_*\le 0$ and $s\ge 0$ then $(t-s)^2\ge t_*^2-2t_*s+s^2$.    Set the parameters as  $t=\Delta u$, $t_*=-\alpha_k |\nabla u|^2 (u+\epsilon)^{-1} - \beta_k |x|^{b} u^q$ and $s=(\alpha_{k+1}+1) (u+\epsilon)^{-1}  |\nabla u|^2 $ to get the following lower bound on the square term that appears in (\ref{Ialpha}),
\begin{eqnarray}\label{square}
  && \left(  \Delta u -(\alpha_{k+1}+1)(u+\epsilon)^{-1}   |\nabla u|^2  \right)^2 \ge \left(\alpha_k |\nabla u|^2 (u+\epsilon)^{-1}+ \beta_k |x|^{b} u^q\right)^2 \\&& \nonumber+2\left(\alpha_k |\nabla u|^2 (u+\epsilon)^{-1}  +\beta_k |x|^{b} u^q\right) (\alpha_{k+1}+1) (u+\epsilon)^{-1}  |\nabla u|^2 \\&&\nonumber +(\alpha_{k+1}+1)^2 (u+\epsilon)^{-2} |\nabla u|^4.
 \end{eqnarray}
 Substitute (\ref{uwk+1}) into the term $ - |\nabla u|^2 (u+\epsilon)^{-2} \Delta u$ that appears in \eqref{Ialpha} to eliminate the Laplacian operator.    Then, apply inequality (\ref{square}) to simplify \eqref{Ialpha} as
\begin{eqnarray*}
&&\frac{I}{\alpha_{k+1}}-2 (u+\epsilon)^{-1} \nabla u\cdot\nabla w_{k+1} \ge \frac{2}{n}(u+\epsilon)^{-1} \{ (\alpha_{k+1}+\alpha_k+1)^2 |\nabla u|^4 (u+\epsilon)^{-2} \\&&+\beta_k ^2 |x|^{2b} u^{2q} + 2\beta_k(\alpha_{k+1}+\alpha_k+1)|x|^b u^{q} (u+\epsilon)^{-1}    |\nabla u|^2    \} - w_{k+1} (u+\epsilon)^{-2} |\nabla u|^2 \\&& -\alpha_{k+1} (2\alpha_{k+1}+1)  |\nabla u|^4 (u+\epsilon)^{-3}
+ \beta_{k+1} |x|^b u^q  (u+\epsilon)^{-2} |\nabla u|^2 + T_{k}.
 \end{eqnarray*}
Collecting similar terms and using the value of  $ T_{k}$,  we end up with
\begin{eqnarray*}
&&\frac{I}{\alpha_{k+1}}-2 (u+\epsilon)^{-1} \nabla u\cdot\nabla w_{k+1} +w_{k+1} (u+\epsilon)^{-2} |\nabla u|^2 \\&& \ge \frac{2}{n} \beta_k^2|x|^{2b} u^{2q} (u+\epsilon)^{-1}+
I^{(2)}_{\alpha_k} |\nabla u|^4 (u+\epsilon)^{-3} + S_{\epsilon,\alpha_k,\beta_k}   |\nabla u|^2 u^{q-2} |x|^b \\&&-2 \beta_{k+1}b |x|^{b-2} (u+\epsilon)^{-1} u^{q} \nabla u\cdot x,
 \end{eqnarray*}
where
\begin{eqnarray*}
I^{(2)}_{\alpha_k}&:=& \frac{2}{n}(\alpha_{k+1}+\alpha_k+1)^2 -2\alpha_{k+1}(\alpha_{k+1}+1) +\alpha_{k+1},\\
S_{\epsilon,\alpha_k,\beta_k} &:=&  \frac{4}{n} \beta_k (\alpha_{k+1}+\alpha_k+1) \frac{u^2}{(u+\epsilon)^2}+ \beta_{k+1} \frac{u^2}{(u+\epsilon)^2} -2\beta_{k+1}q \frac{u}{u+\epsilon}.
 \end{eqnarray*}
Therefore,  the following lower bound for $I$ holds,
\begin{eqnarray}\label{lowerI}
I &\ge& 2 \alpha_{k+1} (u+\epsilon)^{-1} \nabla u\cdot\nabla w_{k+1} 
\\&& \nonumber - \alpha_{k+1} w_{k+1} (u+\epsilon)^{-2} |\nabla u|^2 
\\&& \nonumber + \frac{2}{n}  \alpha_{k+1} \beta_k^2|x|^{2b} u^{2q} (u+\epsilon)^{-1}  + I_{\alpha_k} |\nabla u|^4 (u+\epsilon)^{-3} 
\\&& \nonumber + S_{\epsilon,\alpha_k,\beta_k}   |\nabla u|^2 u^{q-2} |x|^b  -2 \beta_{k+1}b |x|^{b-2} (u+\epsilon)^{-1} u^{q} \nabla u\cdot x.
 \end{eqnarray}
 Finally,  applying this lower bound  for $I$ and the lower bound given for $J$ in \eqref{J}, from (\ref{wk+1}) we get
\begin{eqnarray*}
&&\Delta w_{k+1}-2 \alpha_{k+1} (u+\epsilon)^{-1} \nabla u\cdot\nabla w_{k+1}+ \alpha_{k+1} (u+\epsilon)^{-2} |\nabla u|^2 w_{k+1} -\beta_{k+1} q u^{q-1} |x|^b w_{k+1}\\
&&  \ge |x|^a u^p  \left( 1- q \beta_{k+1}^2  + \frac{2}{n} \alpha_{k+1} \beta_k^2 \frac{u}{u+\epsilon}\right) + \alpha_{k+1} I^{(2)}_{\alpha_k} |\nabla u|^4 (u+\epsilon)^{-3} \\&&+ \left(  \alpha_{k+1} S_{\epsilon,\alpha_k,\beta_k} + \left(q(q-1)-\alpha_{k+1} q \frac{u}{u+\epsilon} \right)\beta_{k+1}\right)|\nabla u|^2 u^{q-2} |x|^b
\\&& +2 b\beta_{k+1}   \left(q- \alpha_{k+1} \frac{u}{u+\epsilon}\right) |x|^{b-2} u^{q-1} \nabla u\cdot x + b\beta_{k+1} (n+b-2) |x|^{b-2} u^q.
 \end{eqnarray*}
Completing the square finishes the proof.

    \hfill $\Box$


\section{Proof of Theorem \ref{mainres} via Iteration Arguments}\label{secapp}

To apply the iteration argument, we need to develop a maximum principle argument for the following equation
\begin{equation}\label{win}
\Delta w-2 \alpha (u+\epsilon)^{-1} \nabla u\cdot\nabla w +\alpha w (u+\epsilon)^{-2} |\nabla u|^2 - \frac{\beta(p+1)}{2} |x|^{\frac{a}{2}} u^{\frac{p-1}{2}} w = f(x)\ge 0 \ \ \ \mathbb{R}^n
   \end{equation}
 that appears in Proposition \ref{propwk}, where $\alpha,\beta$ are positive constants, $u$ is a solution of (\ref{4henon}) and $w,f\in C^{\infty}(\mathbb R^n)$.

    \begin{lemma}\label{boundwt} Suppose that $w$ is a solution of the  differential inequality (\ref{win}) where $u$ is a solution of (\ref{4henon}) and
         \begin{equation}\label{w}
w=\Delta u+\alpha  (u+\epsilon)^{-1} |\nabla u|^2 + \beta  |x|^{\frac{a}{2}}  u^{\frac{p+1}{2}}
   \end{equation}
 for positive constants $\epsilon$, $\alpha$ and $\beta$. Then, assuming that $p+1>2\alpha$ the following holds
  \begin{equation}\label{w+}
   \Delta \tilde w \ge 0  \ \ \ \text{on } \{w\ge 0\} \subset \mathbb{R}^n
   \end{equation}
   where $ \tilde w= (u+\epsilon)^{t} w$ for $t=-\alpha$.
\end{lemma}

\noindent\textbf{Proof:}  Straightforward calculations show that
\begin{eqnarray*}
\Delta \tilde w &=& (u+\epsilon)^{t} \Delta w+ 2t (u+\epsilon)^{t-1} \nabla u\cdot \nabla w 
\\&&+ t (u+\epsilon)^{t-1} w\Delta u+ t(t-1)  w (u+\epsilon)^{t-2} |\nabla u|^2
\end{eqnarray*}
 We  now add and subtract two terms $\frac{\beta(p+1)}{2}  |x|^{\frac{a}{2}}  u^{\frac{p-1}{2}} (u+\epsilon)^{t} w$ and $t w (u+\epsilon)^{t-2} |\nabla u|^2$ to the above identity and collect the similar terms to get
 \begin{eqnarray*}
\Delta \tilde w &=& (u+\epsilon)^{t} \left(  \Delta w+ 2t (u+\epsilon)^{-1} \nabla u\cdot \nabla w - t w (u+\epsilon)^{-2} |\nabla u|^2 - \frac{\beta(p+1)}{2} |x|^{\frac{a}{2}}   u^{\frac{p-1}{2}} w \right) \\&&+  
\frac{\beta(p+1)}{2} |x|^{\frac{a}{2}}  u^{\frac{p-1}{2}} (u+\epsilon)^{t} w + t w (u+\epsilon)^{t-2} |\nabla u|^2 + t (u+\epsilon)^{t-1} w \Delta u 
\\&& +    t (t-1) w (u+\epsilon)^{t-2} |\nabla u|^2.
\end{eqnarray*}
From the fact that   $t=-\alpha$ and $w$ satisfies (\ref{win}) we get
 \begin{eqnarray*}
\Delta \tilde w  \ge   \frac{\beta(p+1)}{2}  |x|^{\frac{a}{2}}  u^{\frac{p-1}{2}} (u+\epsilon)^{t} w +  t (u+\epsilon)^{t-1} w \Delta u +    t^2 w (u+\epsilon)^{t-1}  \frac{|\nabla u|^2}{u+\epsilon}
\end{eqnarray*}
Note that we can eliminate the gradient term using (\ref{w}) that is  $\alpha (u+\epsilon)^{-1} |\nabla u|^2=w-\Delta u- \beta |x|^{\frac{a}{2}}   u^{\frac{p+1}{2}}$. Therefore, after collecting the similar terms we get
\begin{eqnarray*}
\Delta \tilde w & \ge &   \frac{t^2}{\alpha}  w^2 (u+\epsilon)^{t-1}  +  (u+\epsilon)^{t-1} w t \left(1-\frac{t}{\alpha}  \right) \Delta u 
\\&& +  \beta (u+\epsilon)^{t-1} |x|^{\frac{a}{2}}  u^{\frac{p-1}{2}} w \left(  \frac{(p+1)\epsilon}{2}+ u\left(\frac{p+1}{2} -\frac{t^2}{\alpha}\right)  \right)
\\&=:& R_1+R_2+R_3.  
\end{eqnarray*}
We claim that the above three terms $R_1,R_2,R_3$ are nonnegative when $w \ge 0$.   From the fact that $\alpha>0$ one can see that $R_1$ is nonnegative.  From the definition of $t=-\alpha<0$ we have $t(1-\frac{t}{\alpha})=-2\alpha <0$. This together with Proposition \ref{prop0}, that is $\Delta u \le 0$,  confirms that $R_2$ is  nonnegative.   Positivity of $R_3$ is an immediate consequence of the assumptions. In other words, note that $\beta$ is positive and $\frac{p+1}{2}-\frac{t^2}{\alpha}=\frac{p+1}{2}-\alpha$ is also positive based on the assumptions.  This finishes the proof.

    \hfill $\Box$

We now apply Lemma \ref{boundwt} to show that $w$ that is a solution of (\ref{win}) is negative.

    \begin{lemma}\label{boundw} Suppose that $\tilde w$  and $w$ are the same as Lemma \ref{boundwt}.  Let $u$ be a bounded solution of (\ref{4henon}) then $ w \le 0$.
\end{lemma}

\noindent\textbf{Proof:} The methods and ideas that we apply in the proof are motivated by the ones provided by Souplet in \cite{so}.   Multiply (\ref{w+}) with $\tilde w^s_+$ where $s> 0$ is a parameter that will be determined later.  Then, integration by parts over $B_R$  gives us
\begin{equation}\label{eqws}
0 \le \int_{B_R} \Delta \tilde w \tilde w^s_+= - s \int_{B_R} |\nabla \tilde w_+|^2 \tilde w^{s-1}_+ + R^{n-1} \int_{S^{n-1}} \tilde w_r \tilde w^s_+.
\end{equation}
Therefore,
\begin{equation}\label{eqws2}
 \int_{B_R} |\nabla \tilde w_+|^2  \tilde w^{s-1}_+ \le \frac{1}{s(s+1)}  R^{n-1} \int_{S^{n-1}}  (\tilde w^{s+1}_+)_r  =  C(s) R^{n-1} I'(R),
\end{equation}
where $$I(R):= \int_{S^{n-1}} \tilde w^{s+1}_+= \int_{S^{n-1}}  (u+\epsilon)^{-(s+1)\alpha}  w^{s+1}_+ 
$$  
and $C(s)$ is a constant independent from $R$.  Note that $w$ given as $w=\Delta u+\alpha |\nabla u|^2 (u+\epsilon)^{-1} +\beta |x|^{\frac{a}{2}}  u^{\frac{p+1}{2}}$ satisfies $w\ge 0$ if and only if $-\Delta u \le \alpha |\nabla u|^2 (u+\epsilon)^{-1} + \beta |x|^{\frac{a}{2}}  u^{\frac{p+1}{2}}$.  Therefore,
\begin{equation}\label{eq1}
w^{s+1}_+ \le C |\nabla u|^{2(s+1)} (u+\epsilon)^{-(s+1)} +C |x|^{(s+1)a/2} u^{  (s+1)(p+1)/2}
\end{equation}
where $C=C(\alpha,\beta,s)$. Applying this upper bound for $w_+$, we can get an upper bound for $I(R)$ as following.
\begin{eqnarray}\label{I(R)}
I(R) &\le& C  \int_{S^{n-1}} (u+\epsilon)^{-(s+1)(\alpha+1)} |\nabla u|^{2(s+1)} 
\\&& \nonumber + C R^{\frac{s+1}{2}a}   \int_{S^{n-1}} (u+\epsilon)^{-\alpha(s+1)}   u ^{(s+1)(p+1)/2} \\&\le & \nonumber C(\epsilon) \int_{S^{n-1}}  |\nabla u|^{2(s+1)} + C(\epsilon) R^{\frac{s+1}{2}a}    \int_{S^{n-1}} u ^{\frac{s+1}{2}(p+1)} 
\\&=: & \nonumber C(\epsilon) (I_1(R)+I_2(R)).
\end{eqnarray}
  In what follows we show that there is a sequence $R$ such that the two terms $I_1(R)$ and $I_2(R)$ decay to zero, for a fixed $\epsilon$. We start with $I_2(R)$ that includes an  integral of a positive power of $u$ over the sphere.  Due to the boundedness assumption on $u$,  it is straightforward to relate this term to $L^p$ estimates of $u$ over the sphere.   As a matter of fact, if $(s+1)(p+1)>2p$ then from the boundedness of $u $ we have
  \begin{equation}\label{firsts}
 \int_{S^{n-1}}  u ^{\frac{s+1}{2}(p+1)} \le C(n) ||u ||_{L^p(S^{n-1})}^p
 \end{equation}
  and for the case $(s+1)(p+1) \le 2p$  we can perform the H\"{o}lder's inequality to get
    \begin{equation}\label{firstb}
 \int_{S^{n-1}}  u ^{\frac{s+1}{2}(p+1)}  \le C(n,p) ||u ||_{L^p(S^{n-1})}^  {\frac{(p+1)(s+1)}{2}}. \end{equation}
  So, to prove a decay estimate for $I_2(R)$ we need to construct a decay estimate for   $||u ||_{L^p(S^{n-1})}$.  On the other hand, we apply Lemma \ref{sobolev} to get an upper bound for the first term in (\ref{I(R)}) that is $I_1(R)$.  In fact, from Lemma \ref{sobolev} where $i=1$, $\tau=2(s+1)$ and $t=2$ we have
  \begin{eqnarray}\label{du2s}
||D_x u||_{L^{2(s+1)}(S^{n-1})} &\le& C || D_\theta D_x u  ||_{L^2(S^{n-1})}  + C || D_x u  ||_{L^1(S^{n-1})} \\&\le & \nonumber C R ||  D^2_x u  ||_{L^2(S^{n-1})}  + C || D_x u  ||_{L^1(S^{n-1})}
\end{eqnarray}
for $s=\frac{2}{n-3}$. In order to get a decay estimate for $I_1(R)$, we need decay estimates for the two terms in the right-hand side of (\ref{du2s}) which are $||  D^2_x u  ||_{L^2(S^{n-1})}$ and $|| D_x u  ||_{L^1(S^{n-1})}$.

 We now apply the elliptic estimates given in Section \ref{secEst} to provide decay estimates for $||u||_{{L^p(S^{n-1})}}$, $||D_x u||_{{L^1(S^{n-1})}}$ and $ ||D^2_x u||_{L^2(S^{n-1})}$.     To do so we first find appropriate upper bounds for these terms on the ball of radius $R$. Then we use certain comparing measure arguments to construct decay estimates over the sphere.  So,  from Lemma \ref{ellip} when $\tau=2$, we get
  \begin{equation}\label{d2u}
\int_{R/2}^{R}   ||  D_{x}^{2} u ||_{ L^{2}(S^{n-1}) }^2 r^{n-1} dr  \le C \int_{B_{2R}\setminus B_{R/4}}   |\Delta u|^2  +  C R^{-4}    \int_{B_{2R}\setminus B_{R/4}}     u^2 .
\end{equation}
We now apply Corollary \ref{Delta2} and Corollary \ref{u} to get a decay estimate for the right-hand side of (\ref{d2u}) that is
\begin{eqnarray*}
R^{-4}\int_{B_{2R}\setminus B_{R/4}}  u^2&\le& C R^{-4}\int_{B_{2R}\setminus B_{R/4}}  u \le  C R^{-4} R^{n -\frac{a+4}{p-1}} =  C R^{n-\frac{a+4p}{p-1}},\\
       \int_{B_{2R}\setminus B_{R/4}}  |\Delta u|^2 & \le&  C R^{n-\frac{a+4p}{p-1}},
\end{eqnarray*}
where $C$ is independent from $R$. From this and (\ref{d2u}) we obtain the following desired decay estimate on the Hessian operator of $u$
\begin{equation}\label{d2ufinfin}
\int_{R/2}^{R}   ||  D_{x}^{2} u ||_{ L^{2}(S^{n-1}) }^2 r^{n-1} dr  \le C R^{n-\frac{4p+a}{p-1}}.
\end{equation}
Similarly, from Corollary \ref{corgrad} and Lemma \ref{1bound} we have 
\begin{eqnarray}\label{d2ufinfin1}
\int_{R/2}^{R}   ||  D_{x}u ||_{ L^{1}(S^{n-1}) } r^{n-1} dr  &\le & C R^{n-\frac{p+3+a}{p-1} },
\\ \label{d2ufinfin2}
\int_{R/2}^{R}   ||  u ||^p_{ L^{p}(S^{n-1}) } r^{n-1} dr  &\le & C R^{n-\frac{a+4}{p-1}p }.
\end{eqnarray}

Now let's define the following sets.  These sets are meant to facilitate our arguments towards construction of decay estimates for $||u||_{{L^p(S^{n-1})}}$, $||D_x u||_{{L^1(S^{n-1})}}$ and $ ||D^2_x u||_{L^2(S^{n-1})}$. For a large number $M$,  that will be   determined later, define
\begin{eqnarray*}
\label{mu}\Gamma_{1}(R) &:=&\left\{r\ \in(R/2,R); \  ||u||^p_{{L^p(S^{n-1})}}> M R^{-\frac{a+4}{p-1}p}\right\},\\
\label{mdu}\Gamma_{2}(R) &:=&\left\{r\ \in(R/2,R); \  ||D_x u||_{{L^1(S^{n-1})}}  > M R^{-\frac{p+3+a}{p-1} }  \right\},\\
\label{md4u}\Gamma_{3}(R) &:=&\left\{r\ \in(R/2,R); \  ||D^2_x u||^2_{L^2(S^{n-1})}> M R^{-\frac{a+4p}{p-1}} \right\}.
\end{eqnarray*}
We claim that $|\Gamma_{i}(R)|\le R/4$ for $1\le i\le 3$. Using (\ref{d2ufinfin}), we get
\begin{eqnarray*}
 C &\ge &  R^{-n+\frac{a+4p}{p-1}} \int_{R/2}^{R} ||D_{x}^{2} u||^{2}_{L^2(S^{n-1})} r^{n-1} dr
 \\& \ge & N R^{-n+\frac{a+4p}{p-1}} R^{n-1} \int_{R/2}^{R} ||D_{x}^{2} u||^{2}_{L^2(S^{n-1})}  dr
 \\ \\
& \ge&  N M R^{-n+\frac{a+4p}{p-1}} R^{n-1}  \int_{|\Gamma_3(R)|}  R^{-\frac{a+4p}{p-1}} dr
  \\
& \ge&  N M R^{-n+\frac{a+4p}{p-1}}   R^{n-1} |\Gamma_3(R)| R^{-\frac{a+4p}{p-1}}
\\&=& N  M |\Gamma_3(R)| R^{-1} ,
 \end{eqnarray*}
where $N=(1/2)^{n-1}$.  Therefore,  $|\Gamma_3(R)|\le \frac{C}{NM} R$.  Now choosing  $M$ to be large enough that is $M>\frac{4C}{N}$, we get $|\Gamma_{3}(R)|\le R/4$. Similarly, applying (\ref{d2ufinfin1}) and (\ref{d2ufinfin2}), one can show that $|\Gamma_{i}(R)|\le R/4$ for $1\le i\le 2$. Hence,  $|\Gamma_{i}(R)|\le R/4$ for $1\le i\le 3$ while  $\Gamma_{i}(R)\subset (R/2,R)$.  So, we can find a sequence $\tilde R$ such that 
\begin{equation} \label{hatr}
 \tilde R\in (R/2,R)\setminus \bigcup_{i=1}^{i=3}\Gamma_{i}(R)\neq\phi.
\end{equation}
Therefore,  for the sequence $\tilde R$, we obtain
\begin{eqnarray}
\label{upsu} ||u||^p_{{L^p(S^{n-1})}} &\le& M  { R}^{-\frac{a+4}{p-1}p},\\
  \label{upsdu} ||D_x u||_{{L^1(S^{n-1})}}  &\le& M { R}^{-\frac{p+3+a}{p-1} } , \\
\label{upsd2u}  ||D^2_x u||^2_{L^2(S^{n-1})} &\le & M { R}^{-\frac{a+4p}{p-1}}.
\end{eqnarray}
 Substituting  (\ref{upsu}) into (\ref{firsts}) and (\ref{firstb}) we get the following decay estimate on $I_2(R)$ that is
\begin{eqnarray}\label{I2}
I_2(R) &\le& C \chi \{(s+1)(p+1)>2p \} R^{\frac{s+1}{2}a  -\frac{a+4}{p-1}p} 
\\ && \nonumber
+ C \chi \{(s+1)(p+1)\le 2p \} R^{\frac{s+1}{2}a  -\frac{a+4}{p-1}(p+1)\frac{s+1}{2}}
\\ &=& \nonumber C \chi \{(s+1)(p+1)>2p \} R^{-\eta_1}
\\ && \nonumber + C \chi \{(s+1)(p+1)>2p \} R^{-\eta_2},
\end{eqnarray}
where $\chi$ is the characteristic function,  $\eta_1:=a\left( \frac{p}{p-1}-\frac{s+1}{2}  \right)+ \frac{4p}{p-1}>0$ and $\eta_2:=\frac{s+1}{p+1} (ap+2(p+1))>0$. Note that we have used the fact that $\frac{p}{p-1}-\frac{s+1}{2}>0$ because $0< s=\frac{2}{n-3}\le 1 $ when $n \ge 5$.   On the other hand,  substituting  (\ref{upsdu}) and (\ref{upsd2u}) into the Sobolev embedding (\ref{du2s}) we get
\begin{equation}\label{du4fin}
||D_x u||_{L^{2(s+1)}(S^{n-1})} \le  C { R}^{1-\frac{a+4p}{p-1}}+ C { R}^{-\frac{p+3+a}{p-1}} = 2 C  { R}^{-\frac{p+3+a}{p-1}}.
\end{equation}
From this and the definition of $I_1(R)$ we end up with the following decay estimate on $I_1(R)$ that is
\begin{equation}\label{I1}
I_1(R)= \int_{S^{n-1}}  |\nabla u|^{2(s+1)} \le C  { R}^{-\frac{2(p+3+a)(s+1)}{p-1}}= C R^{-\eta_3},
\end{equation}
where $\eta_3:= \frac{2(p+3+a)(s+1)}{p-1}>0$.  Finally from (\ref{I1}) and (\ref{I2}) we observe that $$I(R) \le C R^{-\eta} \ \ \ \text{for all\  \  } R>1,$$
where $\eta:=\min\{\eta_1,\eta_2,\eta_3\}>0$. So,
$I(R)\to 0 $ as $R\to\infty$. Note that as $R\to\infty$ then $\tilde R\to \infty$. Since $I(R)$ is a positive function and converges to zero, there is a sequence such that the functional $I'(R)$ is nonpositive.  Therefore, (\ref{eqws2}) yields 
\begin{equation}\label{eqws33}
 \int_{B_R} |\nabla \tilde w_+|^2  \tilde w^{s-1}_+ \le 0.
\end{equation}
Hence, $\tilde w_+$ has to be a constant. From continuity of $\tilde w$,  we have $\tilde w\equiv C$.  Note that  the constant $C$ cannot be strictly positive. So, $\tilde w_+=0$ and therefore $w_+=0$. This finishes the proof.

\hfill $\Box$

Note that Lemma \ref{boundwt} and lemma \ref{boundw} imply an iteration argument  for the following sequence of functions when $k\ge -1$
\begin{equation}\label{wkkk}
w_k=\Delta u+\alpha_k  (u+\epsilon)^{-1} |\nabla u|^2 + \beta_k  |x|^{\frac{a}{2}}  u^{\frac{p+1}{2}}
   \end{equation}
 as long as the right-hand side of (\ref{wk+1}) stays nonnegative.  For the rest of this section, we construct sequences $\{\alpha_k\}_{k=-1}$ and $\{\beta_k\}_{k=-1}$ such that the right-hand side of (\ref{wk+1}) is nonnegative.  

\subsection{Constructing sequences $\alpha_k$ and $\beta_k$}
In this part, we define sequences $\alpha_k$ and $\beta_k$ needed for the iteration argument.    

\begin{lemma}\label{alphak} Suppose $\alpha_0=0$ and define
 \begin{equation}
\alpha_{k+1}:= \frac{ 4(\alpha_k+1)-n+\sqrt{n(16 \alpha_k^2+24 \alpha_k +n +8)} }{4(n-1)}.
 \end{equation}
Then $(\alpha_k)_k$ is a positive, bounded and increasing sequence that converges to $\alpha:=\frac{2}{n-4}$ provided $n>4$ and $p>1$. Moreover, for this choice of  $(\alpha_k)_k$, one of the sequences of coefficients defined in Proposition \ref{propwk} is zero, i.e. $ I^{(2)}_{ \alpha_k}=0$.
\end{lemma}

\noindent\textbf{Proof:} It is straightforward to show that for any $k\ge 0$ sequences $\alpha_k>0$. Also, direct calculations show that $\alpha_k\to\alpha:=\frac{2}{n-4}$ provided $\alpha_k$ is convergent. Note that $\alpha_1=\frac{4-n+\sqrt{n^2+8n}}{4n-4}<\frac{2}{n-4}$ and by induction one can see that $\alpha_k \le \alpha$ for all $k\ge 0$. In what follows we show that $\alpha_k$ is an increasing sequence. For any $k$ the difference of $\alpha_k$ and $\alpha_{k+1}$ is the following
 \begin{eqnarray*}
 \alpha_{k+1}-  \alpha_{k} &=&  \frac{ \sqrt{n(16 \alpha_k^2+24 \alpha_k +n +8)} - \left( (n-4)+4a_k(n-2)\right) }{4(n-1)} \\&=& \frac{8(n-1)(n-4)(2\alpha_k+1)}{S_{n,k}} \left(  \frac{2}{n-4} -\alpha_k \right)
 \end{eqnarray*}
where $S_{n,k}=\sqrt{n(16 \alpha_k^2+24 \alpha_k +n +8)} + (n-4)+4a_k(n-2) >0 $.   Therefore, from the fact that $\alpha_k\le \alpha=\frac{2}{n-4}$, we get the desired result.

    \hfill $\Box$

Similarly, in what follows we provide an explicit formula for the sequence $\beta_k$.

\begin{lemma}\label{betak} Suppose $\beta_0=\sqrt{\frac{2}{p+1}}$ and define
 \begin{equation}
\beta_{k+1}:= \sqrt{ \frac{2}{p+1} +\frac{4}{(p+1)n} \alpha_{k} \beta_k^2},
 \end{equation}
where $(\alpha_k)_k$ is as in Lemma \ref{alphak}.  Then $(\beta_k)_k$ is a positive, bounded and increasing sequence that converges to $\beta:=\sqrt\frac{2}{(p+1)-c_n}$ where $c_n=\frac{8}{n(n-4)}$ provided $n >4$ and $p>1$. Moreover, for this choice of  $(\alpha_k)_k$ and $(\beta_k)_k$,  one of the sequences of coefficients defined in Proposition \ref{propwk} is strictly positive, i.e.  $I^{(1)}_{0,\alpha_k,\beta_k} > 0$.
\end{lemma}

\noindent\textbf{Proof:}
The sequence $(\beta_k)_k$ for all $k\ge 0$ is positive. Note that boundedness of the sequence $(\alpha_k)_k$ forces the boundedness of the  $(\beta_k)_k$ meaning that $\beta_{k+1} \le \sqrt{\frac{2}{p+1} +\frac{4\alpha}{(p+1)n} \beta_{k}^2}  $ for any $k$.   By straightforward calculations we get
  \begin{eqnarray*}
\beta_{k+1}^2 \le \frac{2}{p+1} \sum_{i=0}^{k+1}  \left(   \frac{4\alpha}{n(p+1)}  \right)^i .
 \end{eqnarray*}
Note that $ \frac{4\alpha}{n(p+1)}=   \frac{8}{n(n-4)(p+1)}<1$ provided $n>4$ and $p>1$. Therefore, $\sum_{i=0}^{\infty}  \left(   \frac{4\alpha}{n(p+1)}  \right)^i <\infty$.   This proves the boundedness of  $(\beta_k)_k$.

Since $(\alpha_k)_{k=0}$ is an increasing sequence,  the sequence $(\beta_k)_{k=0}$ will  be nondecreasing by induction.  Note that $\beta_1=\beta_0$ and $\beta_2=\sqrt{\frac{2}{p+1}+\frac{8}{(p+1)^2n}
\frac{4-n+\sqrt{n^2+8n}}{4n-4}} > \beta_1=\sqrt{\frac{2}{p+1}}$.  Suppose that $\beta_{k-1}\le \beta_k$ for a certain index $k \ge 2$ then we apply the fact that  $\alpha_{k} \ge \alpha_{k-1}$ to show $\beta_{k}\le \beta_{k+1}$. This can be found as a consequence of the following
\begin{eqnarray*}
 \beta_{k+1}-  \beta_{k} &=&  \frac{\beta^2_{k+1}-  \beta^2_{k}}{ \beta_{k+1}+ \beta_{k} } =  \frac{4}{(p+1)n(\beta_{k+1}+ \beta_{k})} (\beta_k^2\alpha_{k}  - \beta_{k-1}^2 \alpha_{k-1} )
 \\&\ge & \frac{4  \alpha_{k-1} (\beta_k  + \beta_{k-1})    }{(p+1)n (\beta_{k+1}+ \beta_{k})} (\beta_k  - \beta_{k-1}).
 \end{eqnarray*}
 So,   $(\beta_k)_k$ is convergent and converges to $\beta:=\sqrt\frac{2n(n-4)}{(p+1)(n-4)n-8}$. Note that $(p+1)n(n-4)>8$ for $p>1$ and $n> 4$. Therefore, $\beta$ is well-defined.

    \hfill $\Box$

Note that based on the definition of the sequences $\{\alpha_k\}_{k=-1}$ and $\{\beta_k\}_{k=-1}$ we concluded that $I^{(1)}_{0,\alpha_k,\beta_k} > 0$ and $ I^{(2)}_{ \alpha_k}=0$.    In the next two lemmata we investigate the positivity of  sequences  $I^{(3)}_{\epsilon,\alpha_{k},\beta_{k}}$ and $I^{(4)}_{a,\epsilon,\alpha_{k},\beta_{k}}$  appeared in (\ref{wk+1}) in Proposition \ref{propwk}.

\begin{lemma}  Set $\epsilon=0$ in $I^{(3)}_{\epsilon,\alpha_{k},\beta_{k}}$  that is defined in Proposition \ref{propwk}. Then,
 \begin{equation}\label{I}
I^{(3)}_{0,\alpha_{k},\beta_{k}}\to I^{(3)}_{0,\alpha,\beta}:=\frac{4}{n}\alpha\beta(2\alpha+1)+\alpha\beta+\beta q(q-3\alpha-1)
   \end{equation}
 as $k\to\infty$. The constant  $I^{(3)}_{0,\alpha,\beta} $  is positive  provided  $p> \frac{n+4}{n-4}$ and $n>4$.
\end{lemma}

\noindent\textbf{Proof:} Note that when $p> \frac{n+4}{n-4}$ and $n>4$, then we have $\frac{p+1}{2} > \frac{n}{n-4}$. As $k\to\infty$, from Lemma \ref{alphak} and Lemma \ref{betak} the sequences $\alpha_k\to\alpha:=\frac{2}{n-4}$ and $\beta_k\to\beta:=\sqrt\frac{2}{(p+1)-c_n}$. Therefore,
\begin{eqnarray*}
 \frac{I^{(3)}_{0,\alpha,\beta}}{\beta} &=& \frac{4}{n}\left(\frac{2}{n-4}\right) \left(\frac{4}{n-4}+1\right)+\frac{2}{n-4}+\frac{p+1}{2}\left(\frac{p-1}{2}-\frac{6}{n-4}\right) \\&=&\left(\frac{p+1}{2}\right)^2-\left(\frac{p+1}{2} \right)\left(\frac{n+2}{n-4}\right) +  \frac{2n}{(n-4)^2}
 \\&=& \left( \frac{p+1}{2}-\frac{n}{n-4}  \right)\left(\frac{p+1}{2}-\frac{2}{n-4} \right)  > 0.
 \end{eqnarray*}

    \hfill $\Box$

Note that $I^{(4)}_{a,\epsilon,\alpha_{k},\beta_{k}}$ appears in (\ref{wk+1}) mainly because of the weight function $|x|^a$.  In other words, we have $I^{(4)}_{0,\epsilon,\alpha_{k},\beta_{k}}=0$, in case of $a=0$.

\begin{lemma} For any $k\ge 0$,
\begin{equation}
I^{(3)}_{0,\alpha_{k},\beta_{k}}  < \beta_{k+1} (\frac{p+1}{2}-\alpha_{k+1})^2,
\end{equation}
provided $p>\frac{n+4}{n-4}$ and $n>4$.     Therefore, for any $a\ge0$ that satisfies the following upper bound
\begin{equation}\label{ak}
 a \le A_k:=\frac{2(n-2) I^{(3)}_{0,\alpha_{k},\beta_{k}}}{   \beta_{k+1} (\frac{p+1}{2}-\alpha_{k+1})^2 -I^{(3)}_{0,\alpha_{k},\beta_{k}}}
 \end{equation}
the sequence $I^{(4)}_{a,0,\alpha_{k},\beta_{k}}$ is positive for any $k$.
\end{lemma}

\noindent\textbf{Proof:}
 Basic calculations show that
 \begin{eqnarray*}
&& \beta_{k+1} (\frac{p+1}{2}-\alpha_{k+1})^2 -I^{(3)}_{0,\alpha_{k},\beta_{k}}\\ &=&
 \beta_{k+1}(\frac{p+1}{2}-\alpha_{k+1})^2 -  \frac{4}{n}\alpha_{k+1} \beta_{k}  (\alpha_{k+1} 
 \\ &&+\alpha_k+1) - \alpha_{k+1} \beta_{k+1}  - \beta_{k+1} \frac{p+1}{2}(\frac{p+1}{2}-3\alpha_{k+1}-1)\\&\ge&
 \beta_{k+1} ( (\frac{p+1}{2}-\alpha_{k+1})^2 -  \frac{4}{n}\alpha_{k+1} (\alpha_{k+1} +\alpha_k+1) - \alpha_{k+1} \\&& - \frac{p+1}{2}(\frac{p+1}{2}-3\alpha_{k+1}-1)   )
 \\&=&  \beta_{k+1} \left ( \frac{n-4}{n} \alpha_{k+1}^2 - \frac{4}{n} \alpha_{k+1}^2 - \frac{4}{n} \alpha_{k+1} + \frac{(p-1)\alpha_{k+1}}{2} +\frac{p+1}{2} \right)
  \end{eqnarray*}
where we have used the fact that $\beta_{k}$ and $\alpha_k$ are increasing sequences in the first and the second inequality respectively.  Therefore,
  \begin{eqnarray*}
&& \beta_{k+1} (\frac{p+1}{2}-\alpha_{k+1})^2 -I^{(3)}_{0,\alpha_{k},\beta_{k}} 
\\&\ge&
  \beta_{k+1} ( \frac{n-4}{n} \alpha_{k+1}^2 + \alpha_{k+1}(   \frac{p-1}{2}  - \frac{4}{n} \alpha_{k+1} ) +\frac{p+1}{2} - \frac{4}{n} \alpha_{k+1} )
  \\&\ge&  \beta_{k+1} \left( \frac{n-4}{n} \alpha_{k+1}^2 + (\alpha_{k+1}+1)(   \frac{p-1}{2}  - \frac{4}{n} \alpha) \right)
  \\&>&0.
  \end{eqnarray*}
Note that in the last inequality we have used the fact that $  \frac{p-1}{2}  - \frac{4}{n} \alpha=   \frac{p-1}{2}  - \frac{4}{n} \frac{2}{n-4}> \frac{4}{(n-4)n} (n-2)>0$, since $p>\frac{n+4}{n-4}$ and $n>4$.

    \hfill $\Box$

\begin{remark}
It would be interesting if a counterpart of  (\ref{newpoint4henon}) could be proved for bounded solutions of the fourth order semilinear equation $\Delta^2 u=f(u)$ under certain assumptions on the arbitrary nonlinearity $f\in C^1(\mathbb R)$.  We expect that such an inequality could be established for some convex nonlinearity $f$.  
\end{remark}

\section{Appendix}

We would like to mention that given the estimates in Lemma \ref{1bound} and Lemma \ref{1boundLap}, one can provide a somewhat simpler proof for Proposition \ref{prop0} as what follows. 
\\
\\
\noindent\textbf{Second Proof for Proposition \ref{prop0}:} From Lemma \ref{1bound}, we have
$\int_{{\mathbb{R}}^n} |x|^{2-n+a} u^p dx < \infty$. Hence we define
the following function 
$$w(x) = \frac{1}{n(n-2)\omega_n}\int_{{\mathbb{R}}^n}
\frac{|y|^a u^p(y)}{|x - y|^{n-2}} dy.$$ 
It is clear that $w(x) \ge
0$ and $\Delta w = - |x|^a u^p$. This implies that for a solution $u$ of (\ref{4henon}), the function  $h(x):= w(x) + \Delta u(x)$ is a 
well defined harmonic function on ${\mathbb{R}}^n$. Thus for any
$x_0 \in {\mathbb {R}}^n$ and any $R > 0$, by the mean value theorem for
harmonic functions, we will have
\begin{eqnarray}\label{hx0}
h(x_0)& := & \int_{\partial B_R(x_0)} h d\sigma
\\
& = &\int_{\partial B_R(x_0)} (w + \Delta u) d\sigma\nonumber\\
& \le & \nonumber\int_{\partial B_R(x_0)}w d\sigma + \int_{\partial B_R(x_0)}
|\Delta u| d\sigma.
\end{eqnarray}

Since $w(x_0) < \infty$, through Tonelli's theorem, we can change
the order of the integrations to see that the first integral on the
right-hand side of (\ref{hx0}) tends to zero as $R \to \infty$ for
all $R$.  To be more precise notice that, up to a constant multiple, the
first integral can be written as
$$ \int_{{\mathbb {R}}^n} \int_{\partial B_R(x_0)} \frac{ d\sigma_x}{|x -
y|^{n-2}} |y|^a u^p(y) dy.$$
Then we use the fact that $\int_{\partial B_R(x_0)} \frac{d\sigma_x}{|x -
y|^{n-2}}  = |y-x_0|^{2-n}$ if $| y -x_0| > R$ and equals to $ R^{2-n}$
if $|y - x_0| < R$. Thus the integral will split into two parts.
Outside part tends to zero as $R \to \infty$ due to the fact that
$w(x_0) < \infty$ while the inside part tends to zero due to the
fact that, by Lemma \ref{1bound}, 
 \begin{eqnarray*}
R^{2-n}\int_{B_R(x_0)} |y|^a
u^p(y) dy &\le& R^{2-n} \int_{B_{R + |x_0|}(0)} |y|^a u^p dy \\ &\le&  C
R^{2-n} (R+|x_0|)^{n -\frac{4p + a}{p-1}}
 \end{eqnarray*}
 tends to zero as $R \to
\infty$.   The second integral will tend to zero for some sequence of
$R$ by Lemma \ref{1boundLap} again. Apply the above inequality to this
sequence to see that $h(x_0) \le 0$. Since $x_0$ is arbitrary, we
have $-\Delta u \ge 0$. 

\hfill $\Box$

\end{document}